\documentclass[12pt]{article}
\usepackage{amssymb,amsthm,amsmath,amsfonts,latexsym,tikz,hyperref}
\usepackage[hmargin=.8in,vmargin=1in]{geometry}
\usepackage{tikz}
\tikzstyle{vertex}=[circle, draw, inner sep=0pt, minimum size=4pt]
\definecolor{darkgreen}{cmyk}{.9,0,.9,.20}

\newtheorem{thm}{Theorem}[section]
\newtheorem{prop}[thm]{Proposition}
\newtheorem{cor}[thm]{Corollary}
\newtheorem{lem}[thm]{Lemma}
\newtheorem{conj}[thm]{Conjecture}
\newtheorem{exa}[thm]{Example}

\newcommand{\da}{\hs{-2pt}\downarrow}
\newcommand{\zba}{\bar{z}}
\newcommand{\fI}{{\mathfrak I}}
\newcommand{\fR}{{\mathfrak R}}
\newcommand{\Ih}{\hat{I}}

\DeclareMathOperator{\Peak}{Peak}
\DeclareMathOperator{\Arg}{Arg}
\DeclareMathOperator{\av}{av}

\newcommand{\ben}{\begin{enumerate}}
\newcommand{\een}{\end{enumerate}}
\newcommand{\ble}{\begin{lem}}
\newcommand{\ele}{\end{lem}}
\newcommand{\bth}{\begin{thm}}
\renewcommand{\eth}{\end{thm}}
\newcommand{\bpr}{\begin{prop}}
\newcommand{\epr}{\end{prop}}
\newcommand{\bco}{\begin{cor}}
\newcommand{\eco}{\end{cor}}
\newcommand{\bcon}{\begin{conj}}
\newcommand{\econ}{\end{conj}}
\newcommand{\bde}{\begin{defn}}
\newcommand{\ede}{\end{defn}}
\newcommand{\bex}{\begin{exa}}
\newcommand{\eex}{\end{exa}}
\newcommand{\barr}{\begin{array}}
\newcommand{\earr}{\end{array}}
\newcommand{\btab}{\begin{tabular}}
\newcommand{\etab}{\end{tabular}}
\newcommand{\beq}{\begin{equation}}
\newcommand{\eeq}{\end{equation}}
\newcommand{\bea}{\begin{eqnarray*}}
\newcommand{\eea}{\end{eqnarray*}}
\newcommand{\bal}{\begin{align*}}
\newcommand{\bce}{\begin{center}}
\newcommand{\ece}{\end{center}}
\newcommand{\bpi}{\begin{picture}}
\newcommand{\epi}{\end{picture}}
\newcommand{\bpp}{\begin{picture}}
\newcommand{\epp}{\end{picture}}
\newcommand{\bfi}{\begin{figure} \begin{center}}
\newcommand{\efi}{\end{center} \end{figure}}
\newcommand{\bprf}{\begin{proof}}
\newcommand{\eprf}{\end{proof}\medskip}

\newcommand{\bsl}{\begin{slide}{}}
\newcommand{\esl}{\end{slide}}
\newcommand{\bfr}{\begin{frame}}
\newcommand{\efr}{\end{frame}}

\newcommand{\hqed}{\hfill \qed}

\newcommand{\eqqed}[1]{$\rule{1ex}{0ex}\hfill{\dil#1}\hfill\qed$}

\newcommand{\ol}{\overline}

\newcommand{\hs}[1]{\hspace{#1}}
\newcommand{\hso}[1]{\hspace{-1pt}}
\newcommand{\vs}[1]{\vspace{#1}}
\newcommand{\qmq}[1]{\quad\mbox{#1}\quad}

\newcommand{\emp}{\emptyset}

\newcommand{\sbe}{\subseteq}

\newcommand{\case}[4]{\left\{\barr{ll}#1&\mbox{#2}\\#3&\mbox{#4}\earr\right.}

\def\<{\langle}
\def\>{\rangle}

\newcommand{\ree}[1]{(\ref{#1})}

\newcommand{\ra}{\rightarrow}

\newcommand{\be}{\beta}
\newcommand{\ga}{\gamma}
\newcommand{\de}{\delta}

\newcommand{\si}{\sigma}
\renewcommand{\th}{\theta}

\newcommand{\bv}{{\bf v}}

\newcommand{\bbC}{{\mathbb C}}

\newcommand{\bbQ}{{\mathbb Q}}
\newcommand{\bbR}{{\mathbb R}}

\newcommand{\bbZ}{{\mathbb Z}}

\newcommand{\cB}{{\cal B}}

\newcommand{\cD}{{\cal D}}

\newcommand{\cF}{{\cal F}}

\newcommand{\fS}{{\mathfrak S}}

\newcommand{\Sb}{\ol{S}}

\DeclareMathOperator{\Des}{Des}

\newcommand{\dil}{\displaystyle}

\begin{document}
\pagestyle{plain}

\title{Descent polynomials
}
\author{
Alexander Diaz-Lopez\\[-5pt]
\small  Department of Mathematics and Statistics, Villanova University,\\[-5pt]
\small 800 Lancaster Avenue, Villanova, PA 19085, USA, {\tt alexander.diaz-lopez@villanova.edu}\\
Pamela E. Harris\\[-5pt]
\small  Mathematics and Statistics Department, Williams College,\\[-5pt]
\small 18 Hoxsey Street, Williamstown, MA 01267, USA, {\tt pamela.e.harris@williams.edu}\\
Erik Insko\\[-5pt]
\small  Department of Mathematics, Florida Gulf Coast University,\\[-5pt]
\small 10501 FGCU Blvd. South, Fort Myers, FL 33965-6565, USA, {\tt einsko@fgcu.edu}\\
Mohamed Omar\\[-5pt]
\small  Department of Mathematics, Harvey Mudd College,\\[-5pt]
\small 301 Platt Boulevard, Claremont, CA 91711-5901, USA, {\tt omar@g.hmc.edu}\\
Bruce E. Sagan\\[-5pt]
\small Department of Mathematics, Michigan State University,\\[-5pt]
\small East Lansing, MI 48824-1027, USA, {\tt sagan@math.msu.edu}
}

\date{\today\\[10pt]
	\begin{flushleft}
	\end{flushleft}}

\maketitle

\begin{abstract}

Let $n$ be a nonnegative integer and $I$ be a finite set of positive integers.  In 1915, MacMahon proved that the number of permutations in the symmetric group $\fS_n$ with descent set $I$ is a polynomial in $n$.  We call this the descent polynomial.  However, basic properties of these polynomials such as a description of their coefficients and roots do not seem to have been studied in the literature.  Much more recently, in 2013, Billey, Burdzy, and Sagan showed that the number of elements of $\fS_n$ with peak set $I$ is a polynomial in $n$ times a certain power of two.  Since then, there have been a flurry of papers investigating properties of this peak polynomial.  The purpose of the present paper is to study the descent polynomial.  We will see that it displays some interesting parallels with its peak relative.  Conjectures and questions for future research are scattered throughout.

\end{abstract}

%
%

\section{Introduction}
\label{sec:int}

For the rest of this paper, we  let $n$ be a nonnegative integer and $I$ be a finite set of positive integers.  
(In Section~\ref{sec:ocg} we will permit $I$ to contain $0$.)
We will also use the notation
\beq
\label{eq:m}
m=\max(I\cup\{0\}),
\eeq
where the presence of zero ensures that $m$ is well defined even when $I$ is empty.  We also use the standard notation $[n]=\{1,2,\dots,n\}$.  More generally, given integers $\ell,n$  we set 
$$[\ell,n]=\{\ell,\ell+1.\dots,n\},$$ 
and similarly for other interval notations.

Denote by $\fS_n$ the symmetric group of permutations $\pi=\pi_1\pi_2\dots\pi_n$ of $[n]$ written in one-line notation.  
Note that we will sometimes insert commas into such sequences for clarity in distinguishing adjacent elements.
The {\em descent set} of $\pi$ is
$$
\Des\pi=\{i \mid \pi_i>\pi_{i+1}\}\sbe [n-1].
$$
Note that a similar definition can be given for any sequence $\pi$ of integers and we will have occasion  to use that level of generality.
Given $I$ and $n>m$, where $m$ is defined by~\ree{eq:m}, we wish to study the set
$$
D(I;n) = \{\pi\in\fS_n \mid \Des\pi = I\},
$$
and  its cardinality
$$
d(I;n)=\# D(I;n).
$$
As an example, if $I=\{1,2\}$ then 
\beq
\label{eq:D12}
D(\{1,2\};n) =\{\pi\in\fS_n \mid \pi_1>\pi_2>\pi_3<\pi_4<\dots<\pi_n\}.
\eeq
It follows that $\pi_3=1$.  Furthermore, one can pick any two integers from $[2,n]$ to be to the left of $\pi_3$. Placing the integers to the left of $\pi_3$ in decreasing order and the remaining ones to the right of $\pi_3$ in increasing order completely determines $\pi$.  Thus
\beq
\label{eq:d12}
d(\{1,2\};n)=\binom{n-1}{2} = \frac{(n-1)(n-2)}{2},
\eeq
which is a polynomial in $n$.  Using the Principle of Inclusion and Exclusion, MacMahon~\cite[Art.\ 157]{mac:ca} proved that this is always the case.
\bth[\cite{mac:ca}]
\label{th:mac}
For any $I$  and all $n>m$ we have that $d(I;n)$ is a polynomial~in~$n$.\hqed
\eth
We   call $d(I;n)$ the {\em descent polynomial} of $I$.
Although this result was proved in 1915, very little work has been done in the intervening years to study these polynomials in more detail.  The purpose of this work is to rectify this oversight.  We also note that since $d(I;n)$ is a polynomial, we can extend its domain of definition to all complex $n$, which will be a useful viewpoint in the sequel.

Another well-studied statistic on permutations is the {\em peak set} defined by
$$
\Peak\pi=\{i \mid \pi_{i-1}<\pi_i>\pi_{i+1}\} \sbe[2,n-1].
$$
It is not true that any set of integers $I\sbe[2,\infty)$ is the peak set of some permutation.  For example, clearly $I$ can not contain two consecutive indices.   Say that  $I$ is {\em admissible} if  there is some permutation $\pi$ with $\Peak\pi=I$.  
For $I$ admissible and $n>m$, consider the set
$$
P(I;n) =  \{\pi\in\fS_n \mid \Peak\pi = I\}.
$$
To illustrate, if $I=\emp$, then
$$
P(\emp;n) = \{\pi\in\fS_n \mid \text{$\pi_1>\dots >\pi_i<\pi_{i+1}<\dots<\pi_n$ for some $1\le i\le n$}\}.
$$
Noting that $\pi_i$ must be $1$, such a permutation is determined by picking some subset of $[2,n]$ to be to the left of $\pi_i$, then arranging those elements in decreasing order, and finally making the rest an increasing sequence to the right of $\pi_i$.  It follows that
$$
\#P(\emp;n) = 2^{n-1},
$$
which is certainly not a polynomial in $n$.
But nearly one hundred years after MacMahon's theorem, Billey, Burdzy, and Sagan~\cite{bbs:pps} proved the following result.
\bth[\cite{bbs:pps}]
\label{th:bbs}
For any admissible $I$ and all $n>m$ we have that
$$
\#P(I;n) = p(I;n) 2^{n-\#I-1},
$$
where $p(I;n)$ is a polynomial in $n$ taking on integer values in the range $(m,\infty)$.\hqed
\eth
As might be expected, $p(I;n)$ is called the {\em peak polynomial} of $I$.  Inspired by this theorem, a number of papers have been written about properties of peak and related polynomials~\cite{bbps:mew,bft:crp,cvdlopz:nps,dlhio:ppp,dlhip:psc,dnpt:psp,kas:mfp}.  It turns out that many of our results about descent polynomials have analogues for peak polynomials.

The rest of this paper is organized as follows.  In the next section we  derive two  recursions for $d(I;n)$ that prove useful in the sequel.  Section~\ref{sec:coe}  is devoted to the study of the coefficients of $d(I;n)$ when expanded in an appropriately centered binomial coefficient basis for the polynomial ring $\bbQ[n]$.  In particular, we  give a combinatorial interpretation for these constants which  permits us to prove a log-concavity result.  We also explore a conjecture that the coefficients of $d(I;n)$ when expanded  in a differently centered basis alternate in sign. In Section \ref{sec:roo}, we  study the roots of the descent polynomial, including those which are complex.  It will be shown that the elements of $I$ are always integral zeros, and progress will be made towards a conjecture about the location of the full set of roots in the complex plane.  Analogues of $d(I;n)$ in Coxeter groups of type $B$ and $D$ are considered in Section~\ref{sec:ocg}.  
We end with a section containing comments and open questions.  There we  present a result that unifies Theorems~\ref{th:mac} and~\ref{th:bbs} using the concept of consecutive pattern avoidance.

%
%

\section{Two recursions}
\label{sec:rec}

In this section we  derive two recursions for  $d(I;n)$.  The first will be useful in a number of ways, for example in determining the degree of $d(I;n)$ and in finding some of its roots.

If $I\neq\emp$, then we let 
$$
I^-= I-\{m\}.
$$
We first express $d(I;n)$ in terms of $d(I^-;n)$ which will permit latter proofs by induction on $m$ or on $\#I$.
\bpr\label{prop:recursion}
If $I\neq\emp$, then
\beq
\label{eq:I^-}
d(I;n)=\binom{n}{m} d(I^-;m)-d(I^-;n).
\eeq
\epr
\bprf
Consider the set $P$ of permutations $\pi\in\fS_n$  that can be written as a concatenation $\pi=\pi'\pi''$ satisfying
\ben
\item $\#\pi'=m$ and $\#\pi''=n-m$, and
\item $\Des\pi'=I^-$ and $\pi''$ is increasing.
\een
We can write $P$ as the disjoint union of those $\pi$ where $\pi_m'>\pi_1''$ and those where the reverse inequality holds.  
So $\#P = d(I;n)+d(I^-;n)$.

On the other hand, the elements of $P$ can be constructed as follows.  Pick $m$ elements of $[n]$ to be in $\pi'$ which can be done in $\binom{n}{m}$ ways.  Arrange those elements to have descent set $I^-$ which can be done in $d(I^-;m)$ ways.  Finally, put the remaining elements in $\pi''$ in increasing order which can only be done in one way.  If follows that $\#P=\binom{n}{m} d(I^-;m)$.  Comparing this with the expression for $\#P$ at the end of the previous paragraph completes the proof.
\eprf

We can use the previous result to provide a new proof of MacMahon's theorem and to also obtain the degree of $d(I;n)$. 
\bth \label{thm:polydegreen}
For all $I$ we have that $d(I;n)$ is a polynomial in $n$ with $\deg d(I;n)=m$.
\eth
\bprf
We prove this by induction on $\#I$.  If $I=\emp$, then $d(I;n)=1$ and the result clearly holds.  For nonempty $I$, we examine~\ree{eq:I^-}.  We have that $\binom{n}{m}$ is a polynomial in $n$ of degree $m$.  Multiplying by the nonzero constant $d(I^-;m)$ does not change this.  And, by induction, $d(I^-;n)$ is a polynomial of lesser degree so that the first term in the difference is dominant.
\eprf
MacMahon also gave an explicit formula for $d(I;n)$ using the Principle of Inclusion and Exclusion.  As a further application of ~\ree{eq:I^-}, we will now rederive this expression.  Before doing so, we set the following notation.
Recall that a {\em composition of $n$} is a sequence of positive integers summing to $n$.
Given a set of positive integers $I=\{i_1<\dots<i_k\}$ and $n>i_k$ it will be convenient to 
let $i_0=0$ and $i_{k+1}=n$.  Now we can form the {\em difference composition}
\begin{equation}\label{delta(J)}
\de(I)=(i_1-i_0,\ i_2-i_1,\ \dots,\  i_{k+1}-i_k).
\end{equation}
To any composition $\de=(\de_1,\dots,\de_k)$ of $n$ we associate the multinomial coefficient
$$
\binom{n}{\de}=\frac{n!}{\de_1!\dots \de_k!}.
$$
Finally, we let $\binom{I}{i}$ be the set  of all $i$-element subsets of $I$.
\bth[\cite{mac:ca}]
\label{PIE}
If $I$ is a  set of positive integers with $\#I=k$, then
\beq
\label{PIEeq}
d(I;n)=\sum_{i\ge0} (-1)^{k-i} \sum_{J\in \binom{I}{i}} \binom{n}{\de(J)}.
\eeq
\eth
\bprf
We proceed by induction on $\#I$. If $I=\emptyset$, then $d(I;n)=1$. In this case the right-hand side of \eqref{PIEeq} is $\binom{n}{\de(\emptyset)}=1$. We assume that the result holds for all sets $I$ with $\#I\leq k$.
Consider $\#I=k+1$ and $m=\max(I)$.  Note that if $\de^-$ is a composition of $m$ then 
$\binom{n}{m}\binom{m}{\de^-}=\binom{n}{\de}$ where $\de$ is $\de^-$  with $n-m$ appended.  Now using this fact, equation~\eqref{eq:I^-}, and the induction hypothesis we have
\begin{align*}
d(I;n)&=\binom{n}{m}\left[\sum_{i\ge0} (-1)^{k-i} \sum_{J\in \binom{I^-}{i}} \binom{m}{\de(J)}\right]
-\sum_{i\ge0} (-1)^{k-i} \sum_{J\in \binom{I^-}{i}} \binom{n}{\de(J)}\\
&=\sum_{i\ge0}(-1)^{k+1-i}\left[\sum_{J\in\binom{I}{i},\;m\in J}\binom{n}{\de(J)}+\sum_{J\in\binom{I}{i},\; m\notin J}\binom{n}{\de(J)}\right]\\
&=\sum_{i\ge0}(-1)^{k+1-i}\sum_{J\in\binom{I}{i}}\binom{n}{\de(J)},
\end{align*}
as desired. 
\eprf

It will be useful to have a recursion  that does not contain any negative terms.  We will see an application of this recursion when we investigate the expansion of $d(I;n)$ in a certain binomial basis.  A similar recursion was used by Diaz-Lopez, Harris, Insko and Omar~\cite{dlhio:ppp} when they proved the peak polynomial positivity conjecture of Billey, Burdzy, and Sagan~\cite{bbs:pps}.  To state our recursion, we  need some notation.

Suppose $I=\{i_1,\dots,i_\ell\}$ where the integers are listed in increasing order.  We define two related sets of positive integers.  Specifically, for $1\le k\le\ell$ we let
$$
I_k=\{i_1,\dots,i_{k-1},i_k-1,\dots,i_\ell-1\}-\{0\},
$$
and
$$
\Ih_k=\{i_1,\dots,i_{k-1},i_{k+1}-1,\dots,i_\ell-1\}.
$$
Note that subtracting $\{0\}$ in $I_k$ is only necessary when $k=1$ and $i_1=1$ so that $I_k$ is still a set of positive integers.  The reason these sets are interesting is that
if one removes $n+1$ from a $\pi\in D(I;n+1)$ then the resulting $\pi'$ has $\Des\pi'=I_k$  or  $\Des\pi'=\Ih_k$ for some $k$. 
 Also note that $n+1$ can only appear at the end of $\pi$ or at a position $i_k$ where $i_k-1\not\in I$.  So define
$$
I'=\{i_k \mid i_k-1\not\in I\}
$$
and $I''=I'-\{1\}$.  Note $I'$ and $I''$ are only different if $1\in I'$.
\bth
\label{thm:rec2}
If $I\neq\emp$, then
$$
d(I;n+1)=d(I;n)+\sum_{i_k\in I''} d(I_k;n) +\sum_{i_k\in I'} d(\Ih_k;n).
$$
\eth
\bprf
We partition $D(I;n+1)$ according to the position of $n+1$.  If $\pi\in D(I;n+1)$ then we let $\pi'$ be the permutation obtained from $\pi$ by deleting $n+1$.  If $\pi_{n+1}=n+1$ then the corresponding $\pi'$ are the elements of $D(I;n)$ which gives the first term in the sum for $d(I;n+1)$.  Now suppose $\pi_{i_k}=n+1$ where $i_k>1$ and $\pi_{i_k-1}>\pi_{i_k+1}$.  Then the possible $i_k$ where this could occur are exactly the elements of $I''$, and the $\pi'$ which result form the set $D(I_k;n)$.  This explains the first summation.  Similarly, suppose  $\pi_{i_k}=n+1$ where either $i_k=1$, or $i_k>1$ and $\pi_{i_k-1}<\pi_{i_k+1}$.  Then the corresponding $\pi'$ are counted by the second sum and we are done.
\eprf

%
%

\section{Coefficients}
\label{sec:coe}

In this section we show that the coefficients of descent polynomials, written in a certain polynomial basis, are nonnegative by providing a combinatorial interpretation for them. Based on a partial result and computer evidence, we then conjecture that these coefficients form a log-concave sequence.   We also make a conjecture that the coefficients in another polynomial basis alternate in sign and prove it in a special case.

The study of coefficients of polynomials has a rich history and many important examples. 
For instance, Ehrhart polynomials \cite{s:drcp} and chromatic polynomials \cite{b:ecpl} can be written in certain polynomial bases using {nonnegative} coefficients. {In 2013} Billey, Burdzy, and Sagan conjectured that peak polynomials {could} be written with non-negative coefficients in a binomial basis \cite{bbs:pps}. This conjecture was {proved} in 2017 by Diaz-Lopez et al. \cite{dlhio:ppp}. We restate their result here and then prove a similar, {but stronger,} result for descent polynomials in Theorem \ref{comb interp}.
\bth [{\cite{dlhio:ppp}}]\label{thm:dlhio} 
For any non-empty admissible set $I$ we have
\[ p(I;n)=b_0(I) \binom{n-m}{0}+b_1(I)\binom{n-m}{1}+\dots+ b_{m-1}(I) \binom{n-m}{m-1},\]
where $b_0(I)=0$ and for $ 1\leq k \leq m-1$ the constant $b_k(I)$ is positive.\hqed
\eth

Before proving our main result of this section, we need a lemma which is of interest in its own right.  Recall the definition 
that for integers $\ell,n$
$$
[\ell,n]=\{\ell,\ell+1,\ell+2,\dots,n\}.
$$
We also use this notation for the sequence $\ell,\ell+1,\dots,n$.  Context should make it clear which interpretation is meant.

\ble
\label{lem:pos}
For any finite set of positive integers $I$ and $n>m$ we have $D(I;n)\neq\emp$.
\ele
\bprf
We induct on $\#I$. If $I=\emp$ then the identity permutation is in $D(I;n)$. If $I\neq\emp$ then by induction there is a permutation $\pi \in \mathfrak{S}_m$ with $\pi\in D(I^-;m)$ where, as usual, $I^-=I-\{m\}$.  It follows that $D(I;n)$ contains the concatenation $\si=\pi' 1 [m+2,n]$ where $\pi'$ is $\pi$ with all its elements increased by one. 
\eprf

We can now state the main result of this section for descent polynomials.

\bth
\label{comb interp}
For any {finite set of positive integers} $I$ we have that
\beq
\label{a_k(I) exp}
d(I;n)=a_0(I) \binom{n-m}{0}+a_1(I)\binom{n-m}{1}+\dots+ a_m(I) \binom{n-m}{m},
\eeq
where $a_0(I)=0$ and for $k\geq 1$ the constant $a_k(I)$ is the number of $\pi\in D(I;2m)$ such that 
\beq
\label{pi cap int}
\{\pi_1,\dots,\pi_m\}\cap [m+1,2m]=[m+1,m+k].
\eeq 
Moreover, $a_k(I)>0$ for $1 \leq k \leq m$.
\eth
\bprf 
By Theorem \ref{thm:polydegreen}, $d(I;n)$ is a polynomial in $n$ of degree $m$, so we can write it uniquely as a linear combination of the polynomial basis 
$$\left\{ \binom{n-m}{0}, \binom{n-m}{1},\dots,\binom{n-m}{m}\right\}. $$
 For ease of notation, given $\pi\in D(I;n)$ we let
$$
\pi[m] = \{\pi_1,\dots,\pi_m\}\cap [m+1,n]. 
$$
Now consider
$$
D_k(I;n)=\{\pi\in D(I;n)\ |\ \#\pi[m]=k\}.
$$
Clearly $D(I;n)$ is the disjoint union of the sets $D_k(I;n)$ for $k\geq 0$.  So to prove the summation formula in~\eqref{a_k(I) exp}, it suffices to demonstrate that $\#D_k(I;n)=a_k(I)\binom{n-m}{k}$.  
We also claim that $D_0(I;n)=\emp$ which forces $a_0(I)=0$.  Indeed, if there is an element $\pi\in D_0(I;n)$ then $\pi[m]=\emptyset$.  This implies that 
$\{\pi_1,\dots,\pi_m\}=[m]$.  Thus $\pi_m\le m$ and $\pi_{m+1}>m$ which contradicts the fact that $m$ is a descent.  

For the rest of the proof we will assume $n\ge 2m$.  This assumption is without loss of generality since if we can show that the polynomials on both sides of equation~\ree{a_k(I) exp} agree for an infinite number of values, then they must agree everywhere.
For $k\geq 1$, consider the elements $\pi \in D_k(I;n)$. There are $\binom{n-m}{k}$ ways to pick the $k$ elements of 
$\pi[m]$.  Furthermore, given any two $k$-element subsets $X$ and $Y$ of $[m+1,n]$, there is an order preserving bijection 
$f:X\ra Y$.   This induces a bijection from the $\pi\in D_k(I;n)$ with $\pi[m]=X$ to the $\si\in D_k(I;n)$ with   $\si[m]=Y$ by applying $f$ to the elements of $\pi[m]$,  leaving the elements in the first $m$ positions from $[m]$ unchanged, and then listing the remaining elements in increasing order.  Note that all the elements of $[m]$ remain unchanged as $f$ is only applied to elements of $[m+1,n]$. This bijection clearly preserves the descent set everywhere except possibly at position $m$.  To see that the descent at $m$ is preserved, note that  $\pi_{m+1}\in[m]$ since the subsequence $\pi_{m+1}\cdots \pi_n$ is increasing
 and there is at least one element of $[m]$ not in  $\{\pi_1,\dots,\pi_m\}$ {because of the assumption} $k\ge1$.  But then in $\si=f(\pi)$ we have $\si_{m+1}=\pi_{m+1}$ since elements of $[m]$ are unchanged.  So if $\pi_m\in[m]$ then $\si_m=\pi_m>\pi_{m+1}=\si_{m+1}$ and if $\pi_m>m$ then $\si_m>m\ge\si_{m+1}$ as desired.

Letting $X=[m+1,m+k]$ we have shown that
 \[\#D_k(I;n)=\#X \cdot \binom{n-m}{k}.\] 
 Furthermore
$k=\#X$ is less than or equal to $m$, which means that the largest interval we need to consider is $[m+1,2m]$ and this is contained in $[m+1,n]$ by our assumption that $n\ge2m$. 
Thus {$\#X=a_k(I)$}  which is clearly a constant independent of $n$.  This completes the proof of the summation formula \eqref{a_k(I) exp}.

To prove the last statement of the theorem, suppose $1\leq k \leq m$. It is enough to show that $D_k(I;2m) \neq \emptyset$. 
By Lemma~\ref{lem:pos} there is $\pi\in D(I^-;m)$.  Thus the concatenation $\si=\pi' [1,k] [m+k+1,2m]$ is in  $D_k(I;2m)$ where $\pi'$ is $\pi$ with every element increased by~$k$.
\eprf

To illustrate  this result, let $I=\{1,2\}$. Then $a_1(I)$ is the number of $\pi=\pi_1\pi_2\pi_3\pi_4 \in D(I;4)$ such that $\{ \pi_1,\pi_2\} \cap [3,4] = [3]$. Similarly,  $a_2(I)$ is the number of $\pi \in D(I;4)$ such that $\{ \pi_1,\pi_2\} \cap [3,4] = [3,4]$. Out of the three elements in $D(I;4)$ one can quickly check that only $\pi=3214$ satisfies the condition for $a_1(I)$, thus $a_1(I)=1$. Similarly, only $\pi=4312$ satisfies the condition for $a_2(I)$, so $a_2(I)=1$. Theorem \ref{comb interp} states that 
\[ d (I;n)= \binom{n-2}{1}  + \binom{n-2}{2}.\] 
By the binomial recursion, this expression agrees with \eqref{eq:d12}.

Many coefficient sequences of combinatorial polynomials have interesting properties, one of which we will investigate in the context of the previous theorem.
A sequence of real numbers $(a_k)=(a_k)_{k\ge0}$ is {\em log-concave} if, for every $k\ge1$, we have $a_{k-1} a_{k+1}\le a_k^2$. Log-concave sequences appear naturally in combinatorics, {algebra, and geometry}; we refer the reader to \cite{s:logconc} and \cite{b:logconc} for important examples and results.  We make the following conjecture about the sequence $(a_k(I))$ which has been verified for any set $I$ with $m\leq 18$.

\bcon 
\label{con:lc}
For any finite set of positive integers $I$,  the sequence $(a_k(I))$ is log-concave. 
\econ

We are able to prove this conjecture for certain $I$,  but first we need a lemma. In it, the sequence $(a_k)$ is said to have a certain property, such as nonnegativity, if all the individual $a_k$ do.  Also, the sequence has {\em no internal zeros} if the elements between any two nonzero elements of the sequence are also nonzero.
\ble
\label{log conc lem}
\ben
\item If $(a_k)$ and $(b_k)$ are log-concave sequences, then so is $(a_k b_k)$.
\item Let $(a_k)$ be a nonnegative log-concave sequence with no internal zeros and let $\ell$ be a positive integer.  Then the sequence 
$(a_k+a_{k+1}+\dots+a_{k+\ell})$ is log-concave.
\een
\ele
\bprf
Statement 1 follows easily from the definition of log-concavity.  For statement 2 note that if we can prove the case $\ell=1$ then the general case will follow by induction since $a_k+\dots+a_{k+\ell}=(a_k+\dots+a_{k+\ell-1})+a_{k+\ell}$.  A nonnegative log-concave sequence $(a_k)$ with no internal zeros satisfies $a_{k+1}/a_{k} \leq a_{k}/a_{k-1}$ for all $k$. 
In particular, if $j\le k$ then $a_{k+1}/a_k \leq a_j/a_{j-1}$  and thus $a_{j-1} a_{k+1}\le a_j a_k$.  So
\begin{align*}
(a_{k-1}+a_k)(a_{k+1}+a_{k+2})
&=a_{k-1}a_{k+1}+a_{k-1}a_{k+2}+a_k a_{k+1} + a_k a_{k+2}\\
&\le a_k^2+a_k a_{k+1}+a_k a_{k+1} + a_{k+1}^2\\
&=(a_k+a_{k+1})^2,
\end{align*}
as desired.
\eprf

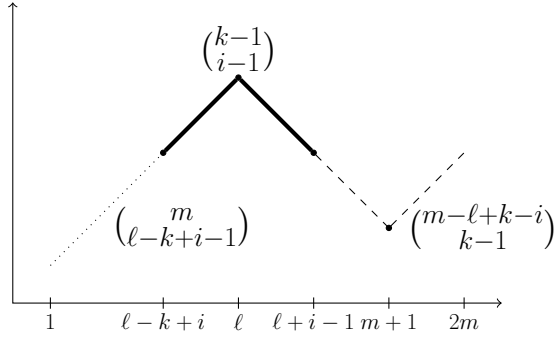
\begin{figure}\label{fig:pi_k}
\begin{center}
	\begin{tikzpicture}[scale=.5, transform shape]
	\node[vertex,fill](a1) at (4,4) {};
	\node[vertex,fill](a2) at (6,6) {};
	\node[vertex,fill](a3) at (8,4) {};
	\node[vertex,fill](a4) at (10,2) {};	
	\draw[ dotted](1,1)--(a1);
	\draw[ ultra thick](a1)--(a2);
	\draw[ ultra thick](a2)--(a3);
	\draw[dashed](a3)--(a4);
	\draw[dashed](a4)--(12,4);
	\draw[<->](0,8)--(0,0)--(13,0);
	\draw[-](1,.15)--(1,-.15);
	\node at (1,-.5) {\Large{$1$}};
	\draw[-](4,.15)--(4,-.15);
	\node at (4,-.5) {\Large{$\ell-k+i$}};
	\draw[-](6,.15)--(6,-.15);
	\node at (6,-.5) {\Large{$\ell$}};
	\draw[-](8,.15)--(8,-.15);
	\node at (8,-.5) {\Large{$\ell+i-1$}};
	\draw[-](10,.15)--(10,-.15);
	\node at (10,-.5) {\Large{$m+1$}};
	\draw[-](12,.15)--(12,-.15);
	\node at (12,-.5) {\Large{$2m$}};	
	\node at (4.5,2) {{\huge{$\binom{m}{\ell-k+i-1}$}}};
	\node at (6,6.75) {{\huge{$\binom{k-1}{i-1}$}}};
	\node at (12.5,2) {{\huge{$\binom{m-\ell+k-i}{k-1}$}}};
\end{tikzpicture}
\caption{The diagram of a $\pi\in D(I;2m)$. The binomial numbers correspond to the possible of ways of choosing each of the three highlighted segments.}
\end{center}
\end{figure}

The next result shows that the sequence $(a_k(I))$ is log-concave in a special case.
\bpr\label{prop:logconc}
Let $\ell\le m$ be positive integers and let $I=\{\ell,\ell+1,\dots,m\}$.  Then $(a_k(I))$ is log-concave.
\epr
\bprf
We first use the combinatorial description of {$a_k(I)$} in Theorem~\ref{comb interp} to derive an explicit formula for this quantity.  Let $\pi\in D(I;2m)$ satisfy equation~\ree{pi cap int}. In Figure \ref{fig:pi_k} we  create a diagram of the permutation $\pi$ by plotting the points $(i,\pi_i)$ and connecting them by, possibly dotted or dashed, segments. Note that the form of $I$ implies that $\pi_1\dots\pi_m$ has a single local maximum at $\pi_\ell$. 
Combining this with~\ree{pi cap int} we see that $\pi_\ell=m+k$ and the elements of $[m+1,m+k]$ are 
$\pi_{\ell-k+i},\pi_{\ell-k+i+1},\dots,\pi_{\ell+i-1}$ for some $i$ with 
$1\le i\le m-\ell+1$.  Now there are $\binom{k-1}{i-1}$ ways of selecting the elements $\pi_{\ell+1},\dots,\pi_{\ell+i-1}$.  Once these elements are put in a decreasing sequence just after $\pi_\ell$, the rest of the elements {of $[m+1,m+k]$} must form an increasing sequence just before $\pi_\ell$.  Next we choose the elements of the increasing sequence $\pi_1,\dots,\pi_{\ell-k+i-1}$ from $[m]$ in $\binom{m}{\ell-k+i-1}$ ways.  The remaining $m-\ell+k-i+1$ elements of $[m]$ must be arranged as the elements $\pi_{\ell+i},\dots,\pi_{m+k}$ with unique local minimum at $\pi_{m+1}$.  So the number of ways to choose $\pi_{m+2},\dots,\pi_{m+k}$ is $\binom{m-\ell+k-i}{k-1}$.  And once these elements are chosen there is only one way to arrange them and the remaining elements since they are all in increasing or decreasing  order.  So
$$
a_k(I)=\sum_{i=1}^{m-\ell+1} \binom{k-1}{i-1}\binom{m}{\ell-k+i-1}{\binom{m-\ell+k-i}{k-1}.}
$$
Now for any fixed $c$, the binomial coefficient sequences $\left( \binom{k}{c} \right)_{k\ge0}$ and $\left( \binom{c}{k} \right)_{k\ge0}$
 are well known to be log-concave.  In combination with Lemma \ref{log conc lem}, this shows that the sequence $(a_k(I))$ is log-concave.
\eprf

If we expand $d(I;n)$ in the binomial basis centered at $-1$ then these coefficients also seems to be well behaved.    The following conjecture has been verified for all $I$ with $m\le 12$.
\bcon
\label{-1 basis}
For any  $I$  we have
$$
d(I;n) = \sum_{k=0}^m (-1)^{m-k} c_k(I) \binom{n+1}{k},
$$
where $c_k(I)$ is a nonnegative integer for all $0\le k\le m$.
\econ

We are able to prove this conjecture for $c_0(I)$.  To do so, we need a couple of lemmas.  Recall that since $d(I;n)$ is a polynomial in $n$, it is defined for all complex numbers.
\ble\label{lem:d at 0}
We have
	\[ d(I;0)=(-1)^{ \#I}.\]
\ele
\bprf
We proceed by induction on $\#I$. The result is clear when $I=\emptyset$ as $d(\emptyset;n)=1$. Consider any set $I$ with $\#I\geq 1$, then by Proposition \ref{prop:recursion} and the inductive hypothesis
	\[ d(I;0)=\binom{0}{m}d(I^-;m)-d(I^-;0) = 0- (-1)^{\#I^-}=(-1)^{\#I},\]
which is what we wished to prove.
\eprf

Keeping the notation of Conjecture~\ref{-1 basis}, we note that 
\beq
\label{c0}
d(I;-1)=(-1)^m c_0(I).
\eeq
This is why our next result will be useful.
\bpr\label{prop:ineqonN}
For any  $I$ and any $n\geq m+2$ we have
\[  d(I;n) \geq |d(I;-1)|.\]
\epr
\bprf
Note that $d(I;n)$  is an increasing function  of $n$ for integral $n>m$ because any permutation $\pi\in D(I;n)$ can be extended to one in $D(I;n+1)$ by merely appending $n+1$.  So it suffices to prove the result when $n=m+2$.

We proceed by induction on $m$. If $m=0$ then $I=\emptyset $ and $d(I;n)=1$ and the result follows.
For the induction step, we first note that by Lemmas \ref{lem:pos} and \ref{lem:d at 0} 
$$
|d(I;0)|=1\leq d(I;m+1).
$$
We now apply Theorem~\ref{thm:rec2}, keeping the notation therein, as well as induction and the previous displayed equation to obtain
	\begin{align*}
		 d(I;m+2)&=d(I;m+1)+\displaystyle\sum_{i_k\in I''}d(I_k;m+1)+\displaystyle\sum_{{i_k\in I'}}d(\Ih_k;m+1)\\
		 	&\geq d(I;m+1)+\displaystyle\sum_{i_k\in I''}|d(I_k;-1)|+\displaystyle\sum_{{i_k\in I'}}|d(\Ih_k;-1)|\\
			&\geq  |d(I;0)|+\displaystyle\sum_{i_k\in I''}|d(I_k;-1)|+\displaystyle\sum_{{i_k\in I'}}|d(\Ih_k;-1)|\\
			&\geq \left|d(I;0)-\displaystyle\sum_{i_k\in I''}d(I_k;-1)-\displaystyle\sum_{{i_k\in I'}}d(\Ih_k;-1)\right|\\
			&= |d(I;-1)|,
	\end{align*}
as desired
\eprf

\bpr
\label{prop:c_0}
For any $I$ we have $c_0(I)\ge0$.
\epr
\bprf
By equation~\ree{c0}, it suffices to show that the sign of $d(I;-1)$ is $(-1)^m$.  We will proceed by induction on $\#I$.  As usual, the case $I=\emp$ is trivial.  For $I\neq\emp$,  applying recursion~\ree{eq:I^-} yields
\begin{equation}\label{eq:dat-1} d(I;-1)=\binom{-1}{m}d(I^-;m)-d(I^-;-1)=(-1)^m d(I^-;m)-d(I^-;-1).
\end{equation}
By Lemma~\ref{lem:pos} we have $d(I^-;m)>0$.
And by induction, the sign of $d(I^-;-1)$ is $(-1)^{m^-}$ where $m^-=\max(I^-\cup\{0\})$. 
So if $m$ and $m^-$ have opposite parity, then  the result follows from \eqref{eq:dat-1}. If  they have the same parity, then $m\geq m^-+2$. Applying Proposition \ref{prop:ineqonN} to $I^-$ we get  $d(I^-;m) \geq |d(I^-;-1)|$.  So, using equation~\eqref{eq:dat-1} again, the sign of $d(I;-1)$ is $(-1)^m$ in this case as well.
\eprf

%
%

\section{Roots}
\label{sec:roo}

We defined $d(I;n)$ only for $n>m$ because we wished to count a nonempty set of permutations.  However, by Theorem~\ref{th:mac}, $d(I;n)$ is a polynomial in $n$ so we can extend the definition to $d(I;z)$ for any complex number $z$.  In this context, it makes sense to talk about the roots of $d(I;z)$ and we study them in this section.  We start by showing that elements of $I$ are roots of $d(I;z)$, a result analogous to one for peak polynomials \cite{bft:crp}.

\bth
\label{roots}
If $I$ is a set of positive integers and $i\in I$ then $d(I;i)=0$. 
\eth
	\begin{proof}
		We induct on $\#I$  using the recursion~\ree{eq:I^-}.  The result is vacuously true when $I$ is empty.  If $i\in I^{-}$ then, by the induction hypothesis, $d(I^{-};i)=0$.  Also $\binom{i}{m}=0$ since $i<m$.  Substituting these values into~\ree{eq:I^-} shows that $d(I;i)=0$.  The only other case is $i=m$.  But then, using equation~\ree{eq:I^-} again, we have that
$$d(I;m)=\binom{m}{m}d(I^{-};m)-d(I^{-};m)=0,$$
as desired.
	\end{proof}

Now that we have established that the elements of $I$ are themselves roots of $d(I;z)$, the remainder of this section focuses on understanding the remaining roots of this polynomial lying in the complex plane.  Throughout we denote by 
$|z_0|$, $\fR(z_0)$ and $\fI(z_0)$ the norm, real and imaginary parts, respectively, of the complex number $z_0$.  

We begin by commenting on the analogous problem for peak polynomials.  Billey, Fahrbach and Talmage \cite{bft:crp} extensively studied the roots of peak polynomials.  Their observations led to the following conjecture regarding the position of the roots in the complex plane.
\bcon[\cite{bft:crp}]
\label{con:peakroots}
For any admissible $I$ and $z_0 \in\bbC$ which is a root of $p(I;z)$, we have
\ben
\item  $|z_0|\le m$, and
\item $\fR(z_0)\ge-3$.
\een
\econ

In fact, in Section 2 of their paper, Billey, Fahrbach and Talmage establish that Theorem~\ref{thm:dlhio} for peak polynomials was implied by this conjecture.  They verified Conjecture~\ref{con:peakroots} computationally for all polynomials $p(I;z)$ where $m \leq 15$.  
We have computed the roots of descent polynomials $d(I;z)$ for all sets $I$ with $m \leq 12$ and arrived at a similar, but more restrictive, conjecture. 

\bcon
\label{con:roots}
For any $I$ and $z_0\in\bbC$ which is a root of $d(I;z)$ we have
\ben
\item  $|z_0|\le m$, and
\item $\fR(z_0)\ge-1$.
\een
\econ

We start by establishing that this conjecture holds for $\#I=1$ by ad hoc means.  Although this approach does not seem to generalize, it gives some intuition about why the two bounds hold.
\bth \label{thm:roots I=m}
If $I=\{m\}$ and $d(I;z_0)=0$ then 
\ben
\item  $|z_0|\le m$, and
\item $\fR(z_0)\ge-1$.
\een
\eth
\bprf
Consider the equation 
$$
0=d(I;z)=\binom{z}{m}-1.
$$
First suppose that $|z|>m$.  Then, by the triangle inequality, $|z-k|\ge |z|-k>m-k$ and  it follows that
$$
\left| \binom{z}{m}\right| = \frac{|z|\cdot |z-1|\cdots |z-m+1|}{m!}>1.
$$
So such $z$ can not be a root of $d(I;z)$ and the first statement in the theorem is proved.

Now suppose $\fR(z)<-1$.  Then $|z-k|\ge |\fR(z-k)|>k+1$ and the previous displayed equation still holds.  This finishes the proof of the second statement.
\eprf

We note that one can use similar techniques to show that if $I=\{1,m\}$ then the roots of $d(I;z)$ satisfy the conjecture. But since we were not able to push this method further we will not present the proof.

In order to establish further bounds for $|z_0|$, we introduce some necessary background on bounding the moduli of roots of polynomials. Recall that given a  nonconstant polynomial $f(z)=\sum_{i=0}^dc_iz^i$, 
the maximum modulus of a root of $f(z)$ is bounded above by the Cauchy bound of $f$, denoted $\rho(f)$, which is the unique positive real solution to the equation
\beq
\label{eq:zero}
|c_0|+|c_1|z+\cdots+|c_{d-1}|z^{d-1}=|c_d|z^d,
\eeq
when $f$ is not a monomial, and zero otherwise
\cite[Theorem 8.1.3]{rs:atp}. 

Although the Cauchy bound of $f(z)$ does not yield an explicit bound for the moduli of the roots of $f(z)$ there are many results that provide such upper estimates for the Cauchy bound. For example \cite[Corollary 8.1.8]{rs:atp} 
gives various bounds for $\rho(f)$ including

\beq
\label{eq:simplebound}
\rho(f)<1+\;\;\max_{0\leq i\leq d-1}\left|\frac{c_i}{c_d}\right|,
\eeq
which we will use in the proof of Theorem \ref{thm:roots}.

It is possible to obtain bounds for polynomials expressed in other bases, such as Newton bases, which we define now.  Given a sequence of complex numbers $\xi_1,\xi_2,\ldots,$ the polynomials 
\[P_k(z)=\prod_{i=1}^{k}(z-\xi_i),\] 
$k\ge0$, form a basis for the vector space of all real polynomials  called the Newton basis with respect to the nodes $\xi_1,\xi_2,\ldots$. Furthermore, since $\deg(P_k(z))=k$ then $\{P_0(z),P_1(z),\ldots,P_d(z)\}$ forms a basis for the vector space of real polynomials of degree at most $d$, for any $d$.
\bth[Theorem 8.6.3  in \cite{rs:atp}]\label{thm:Union}
 Let $f(z)=\sum_{k=0}^dc_kP_k(z)$ be a polynomial of degree $d$ where the $P_k$'s form the Newton basis with respect to the nodes $\xi_1,\ldots,\xi_d$. Then $f$ has all of its zeros in the union of the discs
\beq
\label{eq:discs}
\cD_k:=\{z\in\mathbb{C} \mid  |z-\xi_k|\leq\rho\},
\eeq 
where $k=1,\ldots,d$ and $\rho$ is the Cauchy bound of $\sum_{k=0}^d c_kz^k$. \hqed
\eth
Theorem \ref{thm:Union} played an important role in the work of Brown and Erey that improved known bounds for the moduli of the roots of chromatic polynomials for dense graphs \cite{be:nbcpcr}.  
We will use this result to make progress on Conjecture~\ref{con:roots}.   Because of  recursion~\ree{eq:I^-} we consider the Newton bases with respect to the nodes $0,1,2,3,\ldots$, which is
\[
z\da_k = z(z-1) \cdots (z-k+1), 
\]
$k\ge0$.  This is known as the \emph{falling factorial basis}.  
Expanding $d(I;z)$ in terms of this basis and using the previous theorem immediately gives us the following bounds on the roots of $d(I;z)$.
\ble
\label{lem:CauchyBound}
Suppose $d(I;z)=\sum_{k=0}^{{m}} c_k z\da_k.$  Then the roots of $d(I;z)$ lie in the union of the discs
\[
\cD_k=\{z \in\bbC \mid |z-k|\le \rho(I)\}, 
\]
where $k=0,\ldots,m-1$ and $\rho(I)$ is the Cauchy bound of the polynomial $\sum_{k=0}^{{m}} c_k z^k.$ \hqed
\ele

We now present bounds, linear in $m$, for roots of descent polynomials in the special cases when $\#I \leq 2$, and  bounds which appear to be less tight for general $I$.  We begin by revisitng the case when $\# I = 1$.

\bth
\label{thm:I=m}
Let  $I=\{m\}$ and  
$$
\rho_m=\frac{m}{e}\sqrt[m]{me}.
$$
Then the roots of $d(I;z)$ lie in the union of the discs
$$
\cD_k=\{z\in\bbC \mid |z-k|\le \rho_m\},
$$
where $k=0,\dots,m-1$.
\eth
\bprf
By Lemma~\ref{lem:CauchyBound}, it suffices to show that $\rho(I) \leq \rho_m$.  
Since $d(I;z)=\binom{z}{m}-1$ which has the same roots as $z\da_m - m!$, it suffices show that $\rho_m$ is an upper bound for the unique positive real solution to the equation $z^m=m!$.  This solution is $\sqrt[m]{m!}$, and using upper Riemann sums to estimate the function $\ln m!$ from $\int \ln x \ dx$ establishes that $m!\le m^{m+1}/e^{m-1}$.  The result follows. 
\eprf
We can use the previous result to derive somewhat different bounds from those in Theorem~\ref{thm:roots I=m} for the special case $\# I=1$.
\bco
\label{cor:I=mroots}
If  $I=\{m\}$ and $d(I;z_0)=0$ then
\ben
\item $|z_0|\le \rho_m+m-1$,
\item  $\fR(z_0)\ge-\rho_m$, and
\item $|\fI(z_0)| \le \rho_m$.
\een
Furthermore, for all $m\ge1$, we have
$$
\frac{m}{e} < \rho_m\le m.
$$
\eco
\bprf
Assertions 1, 2 and 3 follow immediately from the description of the discs in Theorem~\ref{thm:I=m}.  To obtain the bounds on $\rho_m$, consider the function $f(m)=\sqrt[m]{me}$.  Taking the derivative gives 
$$
f'(m)=\sqrt[m]{me}\cdot \frac{-\ln m}{m^2} \le 0,
$$
for $m\ge1$.  So $f(m)$ is decreasing  on the interval $[1,\infty)$ and thus is bounded above  by $f(1)=e$.  Applying l'H\^opital's Rule shows that 
$\lim_{m\ra\infty} f(m) = 1$
and this limit is a lower bound.  The desired inequalities follow from observing $\rho_m = mf(m)/e$.  
\eprf

We note that close to the imaginary axis this corollary gives a tighter bound  on $|\fI(z_0)|$
than Theorem~\ref{thm:roots I=m} since $\rho_m\le m$, reducing the area being considered in the earlier theorem by roughly half for large $m$.
We now turn to the case $\#I=2$.

\bth
\label{thm:size2}
Let $I=\{\ell,m\}$ with $1 \leq \ell < m$.  Then the roots of $d(I;z)$ lie in the union of the discs
\[
\cD_k = \{z \in \mathbb{C} \mid |z-k| \leq m\},
\]
for $k=0,\ldots,m-1$.  

\eth

\bprf
We established through computation that the result is true for $m \leq 4$ so we assume $m \geq 5$.
By definition, $I^{-}=\{\ell\}$, so by repeatedly applying equation~\ree{eq:I^-} we have
\begin{align*}
d(I;z) &= \binom{z}{m}d(I^{-};m) - d(I^{-};z) \\
&= \binom{z}{m}\left(\binom{m}{\ell}-1\right) - \binom{z}{\ell}+1 \\
&= \frac{1}{m!}\left(\binom{m}{\ell}-1\right) z\da_m- \frac{1}{\ell!} z\da_{\ell} + 1.
\end{align*}
Multiplying the previous equation by $\ell !$ and using Lemma~\ref{lem:CauchyBound}, the roots of $d(I;z)$ are  contained in the union of the discs
\[\cD_k=
\{z \in \mathbb{C} \mid |z-k| \leq \rho \}, \hspace{0.2in} k=0,1,\ldots,m-1,
\]
where $\rho$ is any upper bound on the unique positive real solution to the equation 
\[
\frac{\ell!}{m!} \left(\binom{m}{\ell} - 1 \right)z^m = z^{\ell}+\ell!.
\]
Since $\binom{m}{\ell} - 1\ge \binom{m}{\ell}/2$, replacing the former expression by the latter in the previous displayed equation just increases the unique positive real solution.  Rewriting the result, it suffices to show that $m$ is an upper bound for the positive real solution of
\[
z^{\ell}\left( \frac{1}{2} \frac{1}{(m-\ell)!}z^{m-\ell} - 1 \right)=\ell!.
\]
To do so, observe that  $m^{\ell}>\ell!$ and $\frac{m^{m-\ell}}{(m-\ell)!}\ge m>4$.  So evaluating the left side of the previous equality at $z=m$ gives
\[
m^{\ell} \left( \frac{1}{2} \frac{1}{(m-\ell)!} \cdot m^{m-\ell}-1 \right) > \ell! \cdot \left(\frac{1}{2} \cdot 4 - 1 \right) = \ell! 
\]
and so $m$ must exceed the unique positive real solution.
\eprf

Similar to Corollary~\ref{cor:I=mroots}, we can use Theorem~\ref{thm:size2} to bound the norm, real and imaginary parts of roots of $d(I;z)$ when $\# I = 2$.

\bco
If $\#I = 2$ and $d(I;z_0)=0$ then
\ben
\item $|z_0|\le 2m-1$,
\item  $\fR(z_0)\ge-m$, and
\item $|\fI(z_0)| \le m$.\hqed
\een
\eco

Similar bounds on the roots of $d(I;z)$ can be established when $\# I = 3$ by first repeatedly applying equation~\ree{eq:I^-} to express $d(I;z)$ as a linear combination of the falling factorials, and then applying a strategy like the one  in the proof of Theorem~\ref{thm:size2}.  But applying these techniques as $\# I$ grows becomes increasingly complicated, so it is not clear that this method will be able to produce a linear bound in general.

We now discuss how to find general bounds on the roots of $d(I;z)$ regardless of the size of $I$.  We begin with the following result.

\ble
\label{le:c_i's}
We have
$$
d(I;z) = c_0 + \sum_{k\in I} c_k z\da_k,
$$
where
$$
\frac{1}{k!} \le |c_k| \le 1,
$$
for all $k\in I\cup\{0\}$.
\ele
\bprf
Induct on $\#I$.  We have $d(\emp;n)=1$ which satisfies the lemma.  By induction we can write
$$
d(I^-;z)= c_0^- + \sum_{k\in I^-} c_k^-  z\da_k,
$$
where
$$
\frac{1}{k!} \le |c_k^-| \le 1,
$$
for all $k\in I^-\cup\{0\}$.
Now using  equation~\ree{eq:I^-} we have that 
\begin{align*}
d(I;z)&= \binom{z}{m} d(I^-;m) - d(I^-;z)\\[5pt]
&=\frac{d(I^-;m)}{m!} z\da_m - \left( c_0^- + \sum_{k\in I^-} c_k^-  z\da_k\right)\\[5pt]
&=-c_0^- - \sum_{k\in I^-} c_k^-  z\da_k + c_m z\da_m,
\end{align*}
where $c_m=d(I^-;m)/m!$.  The lemma now follows for $k<m$ from the bounds on the $c_k^-$, and for $k=m$ from the fact that 
$1\le d(I^-;m)\le m!$.
\eprf

The previous lemma permits us to find general bounds for the roots of $d(I;z)$.

\bth
\label{thm:roots}
Let $I$ satisfy $\# I \geq 2$, and let $m^{-}= \max I^{-}$.  Furthermore let 
\beq
\label{eq:min}
\rho=\min\left(m!+1,\ (m! \cdot \# I)^{1/(m-m^-)}\right). 
\eeq
The roots of $d(I;z)$ all lie in the union of the discs
\beq
\label{m!+1}
\cD_k=\{z\in\mathbb{C}\mid\; |z-k|\le \rho\},
\eeq
where $k=0,1,\dots,m-1$.  In particular, if $d(I;z_0)=0$ then
\ben
\item $|z_0|\le \rho+m-1$,
\item $\fR(z_0)\ge -\rho$, and
\item $|\fI(z_0)|\le \rho$.
\een
\eth
\bprf
The bounds on $|z_0|$, $\fR(z_0)$ and $|\fI(z_0)|$ all follow from~\ree{m!+1}.  Define coefficients $c_k$ as in Lemma \ref{le:c_i's}.  To prove~\ree{m!+1} itself, it suffices to show that $\rho$ is an upper bound for the unique positive real solution of
$$
|c_m| z^m = |c_0|+\sum_{k\in I^-} |c_k|  z^k.
$$
 Replacing $|c_m|$ by its smallest possible value and the other $|c_k|$ by their largest possible value will only increase the value of the positive solution.  So, using the bounds on the $c_k$, it suffices to show that $\rho$ is an upper bound on the unique positive real solution of
\beq\label{eq:auxeq}
\frac{1}{m!} z^m = 1 +\sum_{k\in I^-} z^k.
\eeq
 Applying equation~\ree{eq:simplebound} establishes that $\rho \leq m!+1$.
On the other hand, since $z^k \leq z^{m-}$ for all $k \in I^{-}$ and real $z \geq 1$, $\rho$ is bounded above by the unique positive real solution of the equation $z^m/m! = (\# I) \cdot z^{m^{-}}$, which is $(m! \cdot \# I)^{1/(m-m^-)}$.
\eprf

\begin{figure}
	\centering
		\includegraphics{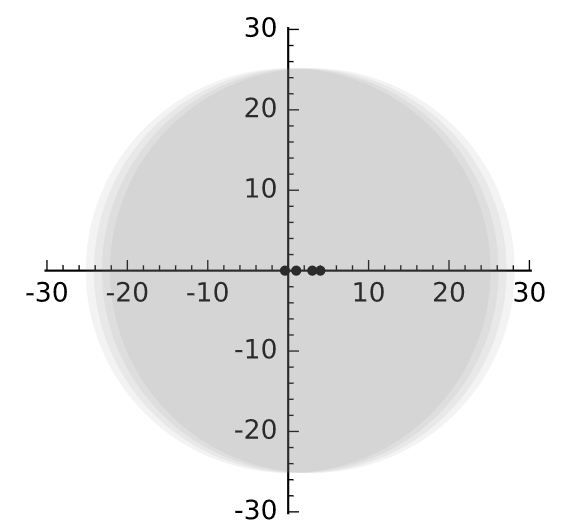} 		\includegraphics{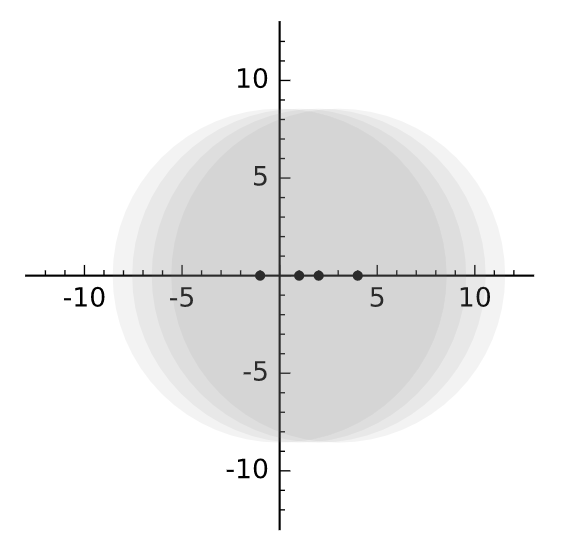}	
	\caption{The roots of descent polynomials for $I=\{1,3,4\}$ and $I=\{1,2,4\}$ are plotted as dots and the corresponding
	bounding discs from Theorem \ref{thm:roots} are shaded in grey.}
\label{fig:boundsm4}
\end{figure}

On the right side of~\ree{eq:min} the first argument achieves the minimum if $m-m^-=1$ since then, using the assumed bound on $\#I$ yields
$$
(m!\cdot \#I)^{1/(m-m^-)}\ge 2 m!>m!+1.
$$
But if $m-m^-\ge2$ then the second argument is smaller since
$$
(m!\cdot \#I)^{1/(m-m^-)}\le (m\cdot m!)^{1/2} < m!+1.
$$
In fact, if $m^-$ is held constant and $m\ra\infty$ then the bound becomes linear.  
An illustration of these two cases is given in Figure~\ref{fig:boundsm4}, where the graph on the left is for $I=\{1,3,4\}$ and the one on the right is for $I=\{1,2,4\}$.

We now use a technique from linear algebra to obtain a different sort of restriction on the roots of $d(I;n)$.  In fact, we will restrict the position of the zeros of any polynomial whose expansion in the falling factorial basis has nonnegative coefficients.  Because of the generality of this result, it will often be less restrictive than Theorem~\ref{thm:roots}.  But along the real axis it will give a linear bound for any $I$ and so it will be an improvement. 

Throughout the remainder of this section, we move freely between a complex number $z=x+iy$ and the vector $\bv=(x,y)\in\bbR^2$.  So if $z=\rho e^{i\th}$ then we call $\th$ an {\em argument} of $\bv$ and write $\arg \bv=\th$. The {\em principle value} of 
$\bv$, denoted $\Arg \bv$, is the argument of $\bv$ satisfying $-\pi<\Arg\bv\le\pi$.  
It will be convenient to let $\Arg (0,0) = \infty$.

Given vectors $\bv_0,\dots,\bv_m$ we say they are {\em nonnegatively linearly independent} if the only linear combination $c_0\bv_0+\dots+c_m\bv_m=(0,0)$ with all the $c_i$ nonnegative is the trivial combination where $c_0=\dots=c_m=0$.   Otherwise the vectors are {\em nonnegatively linearly dependent}.  An {\em open half-plane} consists of all points on one side of a line $L$ through the origin.  The corresponding {\em closed half-plane} is obtained by also including the points on $L$.
The easy backward direction of the following lemma is  a well-known tool in the literature.  But we present a proof for completeness as well as showing that the two statements are actually equivalent.
\ble
\label{nli} 
Vectors $\bv_0,\dots,\bv_m$ are  nonnegatively linearly independent if and only if they all lie in some open half-plane.
\ele
\bprf
If the vectors all lie in an open half-plane then clearly so will any nontrivial nonnegative linear combination.  Since the half-plane is open, such a linear combination can not be zero.

Now suppose the vectors do not lie in any open half-plane.  There are two cases.  If they all lie in a closed half-plane then, since they do not lie in any open half-plane, there must be two of the vectors, say $\bv_0$ and $\bv_1$, such that 
$\bv_0=-c \bv_1$ for some scalar $c>0$.  Thus $\bv_0+c\bv_1=(0,0)$ and the vectors are nonnegatively linearly dependent.

Now suppose that the vectors do not lie in any closed half-plane and consider the vector $\bv_0$.  We will find two other vectors satisfying a nonnegative linear dependence relation with $\bv_0$.   Rotating each of $\bv_0,\dots,\bv_m$ through the  angle 
$-\Arg \bv_0$, we can assume that $\bv_0$ lies along the positive $x$-axis.  Since all the vectors do not lie in the half-plane $x\ge0$ there must be some vector, say $\bv_1$, with $|\Arg\bv_1|>\pi/2$.   Consider the line $L$ through $\bv_1$.  Note that by construction, $\bv_0$ and the negative $x$-axis are on opposite sides of $L$.  And, by the closed half-plane hypothesis again, there must be some $\bv_2$ on the same side of $L$ as the negative $x$-axis but on the opposite side of the $x$-axis from $\bv_1$.  It follows that there is some nonnegative linear combination $a\bv_1+b\bv_2$ which lies on the negative $x$-axis.  So $a\bv_1+b\bv_2=-c\bv_0$ for $c> 0$ which gives the nonnegative linear dependency $c\bv_0+a\bv_1+b\bv_2=(0,0)$.
\eprf

Since the linear dependencies in the previous proof only involve at most three vectors, we have actually proved the following result.
\ble
Vectors $\bv_0,\dots,\bv_m$ are  nonnegatively linearly independent if and only any three of them lie in an open half-plane.\hqed
\ele

To make the connection with roots of polynomials,  let $P_m(z)$ be the vector space of polynomials in a variable $z$ with real coefficients and let $\cB(z)=\{b_0(z),\dots,b_m(z)\}$ be a basis for $P_m(z)$.  Consider the subset of $P_m(z)$ defined by
$$
P_\cB(z)=\left\{f(z)\neq 0  \mid \text{$f(z) =\sum_{k=0}^m c_k b_k(z)$ with $c_k\ge0$ for all $k$}\right\},
$$
where in the above definition $0$ represents the zero polynomial. Translating Lemma~\ref{nli} into this language we immediately have the following result.
\bco
\label{P_B(z)}
The complex number $w$ is not a root of any polynomial in $P_\cB(z)$ if and only if the vectors corresponding to the complex numbers in $\cB(w)$ lie in some open half-plane.\hqed 
\eco
We now specialize to the falling factorial basis $\{z\da_k  | \ k\ge 0\}$.  As usual $\zba$ denotes the complex conjugate of $z$, and if $S$ is a set of complex numbers, then we let
$\Sb=\{\zba \mid z\in S\}$.
\bth
\label{cF}
Let 
$$
\cF(z)=\{z\da_0,\dots,z\da_m\}.
$$
The complex number $w$ is not a root of any polynomial in $P_\cF(z)$ if and only if $w$ is in the region $R=S\cup\Sb$ where
\beq
\label{region}
S=\left\{z\in \bbC \mid \text{$\Arg z \ge 0$ and $\sum_{i=1}^m \Arg(z-i+1)<\pi$}\right\}. 
\eeq

\eth
\bprf
Since the coefficients of polynomials  $f(z)\in P_\cF(z)$ are real, we have $f(w)=0$ if and only if $f(\bar{w})=0$.  So, letting $R$ be the region of $w$ which are not roots of any such $f(z)$, we have $R=S\cup\Sb$ where $S=\{z\in R \mid \Arg z\ge0\}$.  So it suffices to show that $S$ is given as in the statement of the theorem.  Equivalently, by the previous corollary, we must show that
for $z$ with $\Arg z\ge 0$ we have $z\in S$ as defined by equation~\ree{region} if  and only if the elements of $\cF(z)$ lie in an open half-plane.

Suppose first that the sum inequality in~\ree{region} holds for $z$.  Since $z\da_0=1$, we wish to show that for $1\le k \le m$ the complex numbers  $z\da_k$ lie either on the positive $x$-axis or  in the open half-plane above the $x$-axis.  
For then the elements of  $\cF(z)$ will lie in the  open half-plane above the line $y=\epsilon x$ for a sufficiently small negative $\epsilon$.
Since  $\Arg z\ge 0$, we have $\Arg(z-r)\ge 0$ for all reals $r$.  Using this  and the fact that $1\le k\le m$, we have
$$
0 \le \sum_{i=1}^k \Arg(z-i+1) \le \sum_{i=1}^m \Arg(z-i+1)<\pi.
$$
But $z\da_k=\prod_{i=1}^k (z-i+1)$, so the displayed inequalities imply $0\le \Arg(z\da_k)<\pi$ which is what we wished to show.

To complete the proof we must show that if $\sum_{i=1}^m \Arg(z-i+1)\ge\pi$ then the elements of $\cF(z)$ will not all lie in any open half-plane.  From the argument in the preceding paragraph we see that $s_k:=\sum_{i=1}^k \Arg(z-i+1)$ is an increasing function of $k$.  And $s_0=0$.  Thus there must be a nonnegative  integer $\ell$ such that $s_\ell<\pi\le  s_{\ell+1}$.
If $s_{\ell+1}=\pi$ then $z\da_0$ and $z\da_{\ell+1}$ are nonnegatively linearly dependent and we are done by Lemma~\ref{nli}.  If $s_{\ell+1}>\pi$ then we must have $0<\Arg z <\pi$.  It follows that $0< \Arg(z-\ell)<\pi$.
Since $z\da_{\ell+1}=(z-\ell) z\da_\ell$, the previous inequalities force a point on the negative $x$-axis to be a nonnegative linear combination of $z\da_\ell$ and $z\da_{\ell+1}$.  So, together with $z\da_0=1$ we have a nonnegative linear dependency in this case as well.  This concludes the proof of the theorem.
\eprf

Finally, we return to descent polynomials.  If $S$ is any set of complex numbers and $w\in\bbC$ then let
$S+w =\{z+w\mid z\in S\}$.
\bco
\label{co:lin alg}
Let $I$ be a finite set of positive integers.  Then any element of $R+m$ where $R$ is defined as in Theorem~\ref{cF} is not a root of $d(I;z)$.
\eco
\bprf
By Theorem~\ref{comb interp}, we can write
$$
d(I;z)=\sum_{k=0}^m a_k(I) \binom{z-m}{k}=\sum \frac{a_k(I)}{k!} (z-m)\da_k,
$$
where $a_k(I)/k!\ge0$ for all $k$.  So $f(z):=d(I;z+m)\in P_\cF(z)$.  Applying the previous theorem and using the fact that $z\in R+m$ if and only if $z-m\in R$ finishes the proof.
\eprf

\begin{figure}
	\centering
		\includegraphics[height=2in]{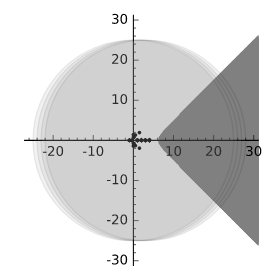} \hspace{.5cm}  \includegraphics[height=2in]{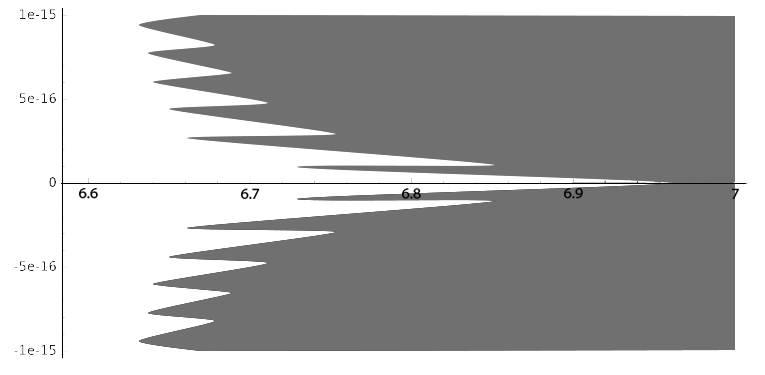}
	\caption{Roots of descent polynomials $d(I;n)$ with $I \subseteq [4]$ plotted inside the
	two bounding regions and close-up view of the region $R+4$ near the real axis.}
	\label{fig:RootBounds411and415cropped}
\end{figure}

Figure \ref{fig:RootBounds411and415cropped} plots all of the roots of descent polynomials corresponding to subsets $I \subseteq [4]$ as small dots,
the worst-case bounds described in Theorem \ref{thm:roots} for such roots are shaded in light grey and  
the dark grey arc is the region $R+4$ where $R$ is as described in Theorem \ref{cF}. 
The image on the right gives a close-up view of the region $R+4$ near the real-axis. 
While in the first image the region $R+4$ looks to be bounded by a curve 
passing through the real-axis near $z=6.65$, it actually passes through the real-axis at $z=7$ 
and then curves back to include complex numbers whose real parts are less than $7$.

\vspace{20pt}

We can use the previous corollary to get our best bound for the size of roots along the positive $x$-axis which holds for general $I$.
\bpr
If $z_0$ is a real root of $d(I;z)$ then $z_0\le 2m-1$.
\epr
\bprf
For a real number $z_0$ we have $\Arg z_0 =0$  if $z_0>0$ and $\Arg z_0=\pi$ if $z_0<0$.  So to be in the region $S$ of equation~\ree{region} we must have $z_0>m-1$.  Applying Corollary~\ref{co:lin alg} we see that if $z_0>2m-1$ then it can not be a zero of $d(I;z)$ and the result follows.
\eprf

%
%

\section{Other Coxeter groups}
\label{sec:ocg}

Recall that for any finite Coxeter system $(W,S)$, the {\em (right) descent set} of $w\in W$ is
\beq
\label{Des w}
\Des w=\{s\in S\ |\ \ell(ws)<\ell(w)\},
\eeq
where $\ell$ is the length function.  
In this section we will consider the Coxeter groups $B_n$ and $D_n$. We will use symbols near the beginning of the Greek alphabet for elements of $B_n$ and $D_n$ to distinguish them from the permutations in $A_{n-1}=\fS_n$.  

We view $B_n$ as the group of signed permutations $\be=\be_1\dots\be_n$ where $\be_i\in\{\pm1,\dots,\pm n\}$ for all $i \in \mathbb{Z}$  
and the sequence $|\be_1|\dots|\be_n|$ is a permutation in $A_{n-1}$, and we view $D_n$ as the subgroup of $B_n$ consisting of all $\be=\be_1\dots\be_n$ where there are  an even number of $\beta_i$ in $\{-1,-2,\dots, -n\}$.
Since $D_n$ is a subgroup of $B_n$, the notation defined below in terms of $B_n$ also applies to $D_n$. We will use the common convention that  $-b$ will be written as $\bar{b}$. 
For example two elements of $B_6$ are $\be=\bar{3} 4 \bar{1} \bar{5} 6 2$ and $\ga= \bar{3}\bar{4}\bar{1}\bar{5} 6 2$, and  the second element is also an element of the subgroup $D_6$, whereas the first is not. 

The simple reflections in $B_n$ are  $S_B=S_A\cup\{s_0\}$ where $s_0=(1,\bar{1})$ and $S_A$ denotes the set of  adjacent  transpositions generating the Coxeter group of type $A_{n-1}$. 
 Identifying reflections and subscripts  as we have done in the symmetric group, we see that for $\be\in B_n$ we have $\Des\be\sbe[n-1]\cup\{0\}$.  
Because of this, it will be convenient to extend permutations in $B_n$ by writing $\be=\be_0\be_1\dots\be_n$
where $\be_0=0$. In this notation, our  previous  examples would be written  $\be=0 \bar{3} 4 \bar{1} \bar{5} 6 2$ and $\ga=0 \bar{3}\bar{4}\bar{1}\bar{5} 6 2$.
Translating definition~\ree{Des w} using our conventions, we see that if $\be=\be_0\be_1\dots\be_n\in B_n$ then
\beq
\label{Des be}
\Des\be = \{i\ge0\ |\ \be_i>\be_{i+1}\} ,
\eeq
where we are using the usual order on the integers for the inequalities.  
To continue our examples in $B_6$, we have $\Des\be=\{0, 2, 3, 5\}$ and $\Des\ga=\{0, 1, 3, 5\}$.

Now given a finite set of nonnegative integers $I$ and $n>m$ where $m$ continues to be defined by   equation~\ree{eq:m},   we let
\beq
\label{D_B def}
D_B(I;n)=\{\be\in B_n\ |\ \Des(\be)=I\} \qmq{and} d_B(I;n)=\#D_B(I;n).
\eeq
We will first derive a recursive formula for $d_B(I;n)$ analogous to the one for $d(I;n)$ in Proposition~\ref{prop:recursion}.

	\bth\label{recursiveBn}
	Let $I$ be a nonempty, finite set of nonnegative integers. Then we have

\vspace{2pt}

\beq 
d_B(I;n)=\binom{n}{m} 2^{n-m} d_B(I^-;m) - d_B(I^-;n). \label{eq:recB}
\eeq
	\eth

	\bprf
Consider the set $P$ of signed permutations $\be \in B_n$ which can be 
written as a concatenation $\be=0\be'\be''$ satisfying
\ben
\item $\#\be'=m$ and $\#\be''=n-m$, and
\item $\Des\be'=I^-$ and $\be''$ is increasing.
\een
We can write $P$ as the disjoint union of those $\be$ where $\be_m'>\be_1''$ 
and those where the reverse inequality holds.  
So $\#P = d_B(I;n)+d_B(I^-;n)$.

On the other hand, the elements of $P$ can be constructed as follows.  
Pick a subset $S$ of $m$ elements  of $[n]$ which can be done in $\binom{n}{m}$ ways.
Form a signed permutation from the elements of $S$  whose descent set is $I^-$ which can be done 
in $d_B(I^-;m)$ ways. Next choose the sign of the $n-m$ elements in $[n] - S$ 
which can be done in $2^{n-m}$ ways. 
Then arrange them in increasing order to form $\be''$ which can only be done in only one way.  
It follows that $\#P=\binom{n}{m}2^{n-m} d_B(I^-;m)$.  Comparing this with the expression for $\#P$ at the end of the previous paragraph completes the proof.
\eprf

Next we prove the type $B$ analogue of Theorem~\ref{PIE}. To state it, we let 
$$I^+=I - \{0\}.$$  
 Also, if $J$ is a set of positive integers then we will let $\de_1(J)$ denote the first component of the composition $\de(J)$.  Note that
 $$
 \delta_1(J)=\case{\min J}{if $J\neq\emp$,}{n}{if $J=\emp$.}
 $$

\bth
\label{PIEB}
If $I$ is a set of nonnegative integers with $\# I^+=k$, then
\beq
d_B(I;n)=\sum_{i\ge0} (-1)^{k-i} \sum_{J\in \binom{I^+}{i}} \binom{n}{\de(J)} \cdot \case{2^{n-\de_1(J)}}{if $0\not\in I$,}{(2^n-2^{n-\de_1(J)})}{if $0\in I$.} \label{eq:Balt}
\eeq
\eth
 
\bprf
We first consider the case where $0 \notin I$ so that $I=I^+$, and proceed by induction on $\#I$. 
If $I=\emptyset$, then  $d_B(I;n)=1$. In this case, the right-hand side of equation \eqref{eq:Balt} also gives 
$\binom{n}{\de(\emptyset)}=1$. 
We assume that the result holds for all sets $I$ not containing $0$ with $\#I\leq k$. 
Consider $\#I=k+1$ and $m=\max(I)$.  Using recursion~\eqref{eq:recB}, and the induction hypothesis we have
\begin{align*}
d_B(I;n)&=\binom{n}{m}2^{n-m} \left[\sum_{i\ge0} (-1)^{k-i} \sum_{J\in \binom{I^-}{i}} \binom{m}{\de(J)}2^{m-\de_1(J)}\right]
-\sum_{i\ge0} (-1)^{k-i} \sum_{J\in \binom{I^-}{i}} \binom{n}{\de(J)}2^{n-\de_1(J)}\\
&=\sum_{i\ge0}(-1)^{k+1-i}\left[\sum_{J\in\binom{I}{i},\;m\in J}\binom{n}{\de(J)} 2^{n-\de_1(J)}+\sum_{J\in\binom{I}{i},\; m\notin J}\binom{n}{\de(J)}2^{n-\de_1(J)}\right]\\
&=\sum_{i\ge0}(-1)^{k+1-i}\sum_{J\in\binom{I}{i}}\binom{n}{\de(J)}2^{n-\de_1(J)}.
\end{align*}
Since $I=I^+$ when $0 \notin I$, this completes the proof for this case.

Next we consider when $0 \in I$.  If $I=\{0\}$ then Theorem~\ref{recursiveBn} shows
$d_B(I;n)=2^n-1$, and the right hand of equation~\eqref{eq:Balt} above gives
$\binom{n}{\delta(\emptyset)}(2^n-2^{n-n})$.  
So equation~\eqref{eq:Balt} holds in this case.
The  induction argument is exactly the same as that of the case when $0 \notin I$, but one replaces $2^{m-\de_1(J)}$ with $2^m-2^{m-\de_1(J)}$ and   $2^{n-\de_1(J)}$ with $2^n-2^{n-\de_1(J)}$.
\eprf

Using Theorems~\ref{PIE} and~\ref{PIEB},
we can also give a simple numerical relationship between the descent formulas in types $A$ and $B$.

\bco
 \label{descentBA}
Let $I$ be a finite set of positive integers and $I_0 = I \cup \{0\}$.  Then

\vs{5pt}

\eqqed{
d_B(I;n)+d_B(I_0;n) = 2^n d(I;n).
}
\eco

Since the right-hand side of equation~\ree{eq:Balt} is well defined for all real numbers $n$, we use it to extend the definition $d_B(I;n)$ to $\bbR$ and talk about its roots.  The proof of the following theorem is similar to that of Theorem~\ref{roots} and so is omitted.

\bth\label{rootsB}
If $I$ is a set of nonnegative integers and $i\in I$ then $d_B(I;i)=0$. \hqed
\eth

The remaining results of this section pertain to the Coxeter group $D_n$.  We continue to use all the conventions for $B_n$ with this subgroup.  In particular, we will use the same definition of  $\Des\be$ as in equation~\ree{Des be}, and the notation $D_D(I;n)$ and $d_D(I;n)$ is defined exactly as in equation~\ree{D_B def} except that $\be$ runs over $D_n$ rather than $B_n$.  Our results in type $D_n$ are very similar to those in type $B_n$ except with some changes imposed by using a different power of two and the intermingling of $d_D$ and $d_B$ in the same formula.

\bth\label{recursiveDn} 
Let $I$ be a nonempty, finite set of nonnegative integers. Then 
\vspace{2pt}
\beq 
d_D(I;n)=\binom{n}{m} 2^{n-m-1} d_B(I^-;m) - d_D(I^-;n). \label{eq:recD}
\eeq
\eth

	\bprf
Consider the set $P$ of signed permutations $\be \in D_n$ satisfying the same two conditions as in the proof of Theorem~\ref{recursiveBn}.
As before,  $\#P = d_D(I;n)+d_D(I^-;n)$.

An alternative construction of the elements of $P$ is as follows.  
Pick $m$ elements from $[n]$ which can be done in $\binom{n}{m}$ ways.  Use those elements to create a type $B$ signed permutation $\beta'$ with descent set $I^-$ which can be done 
in $d_B(I^-;m)$ ways. 
Since a type $D_n$ permutation must have an even number of negative signs, of the remaining $n-m$ elements 
choose the sign of the first $n-m-1$ of them; the sign of the last element in the set of numbers appearing in $\beta''$ is then determined by the number of negative signs  assigned previously.  
Thus choosing the signs of the elements appearing in $\be''$ can be done in $2^{n-m-1}$ ways. 
Now form the unique increasing arrangement of these signed integers to form $\beta''$.
It follows that $\#P=\binom{n}{m}2^{n-m-1} d_B(I^-;m)$ and we are done as in the proof of Theorem~\ref{recursiveBn}.
\eprf

Next we can use Theorem \ref{recursiveDn} to prove a Type $D_n$ analogue of Theorems~\ref{PIE} and \ref{PIEB}.  
  As the proof are similar to those we have seen before, we omit them.

\bth
\label{PIED}
If $I$ is a set of nonnegative integers with $\# I^+=k$, then
$$
d_D(I;n)=
\case{\dil (-1)^k+\sum_{i>0} (-1)^{k-i} \sum_{J\in  \binom{I^+}{i}} \binom{n}{\de(J)} \cdot 2^{n-\de_1(J)-1}}{if $0\not\in I$,}
{\dil (-1)^k(2^{n-1}-1)+\sum_{i>0} (-1)^{k-i} \sum_{J\in \binom{I^+}{i}} \binom{n}{\de(J)} \cdot (2^{n-1}-2^{n-\de_1(J)-1})}{if $0\in I$. \rule{0pt}{30pt}}
$$
forall $n>m$\hqed
\eth

Finally we present the analogues of Corollary~\ref{descentBA}, and Theorem~\ref{rootsB} for type $D_n$.  

 	\bco\label{recursiveDn1}
	Let $I$ be a nonempty set of positive integers and $I_0=I\cup\{0\}$.  Then
	\begin{enumerate}
		\item $d_D(I;n)+d_D(I_0;n) = 2^{n-1} d(I;n)$,  and
		\item  $d_D(I;i)=d_D(I_0;i)=0$  whenever $i\in I^-$. \hqed
	\end{enumerate}
	\eco

%
%

\section{Comments and open questions}
\label{sec:coq}

We end with some comments about our results.  These include avenues for future research and more conjectures.

\medskip

{\bf (1)  Consecutive pattern avoidance.}
One way to unify Theorems~\ref{th:mac} and~\ref{th:bbs} is through the theory of consecutive pattern avoidance.  Call two sequences of integers $a_1 a_2\dots a_k$ and $b_1 b_2\dots b_k$ {\em order isomorphic}  provided $a_i<a_j$ if and only if $b_i<b_j$ for all pairs of indices $1\le i,j\le k$.  Given $\si\in\fS_k$ called the {\em pattern}, we say that $\pi\in\fS_n$ {\em contains a consecutive copy of $\si$ at index $i$} if the factor $\pi_i\pi_{i+1}\dots\pi_{i+k-1}$ is order isomorphic to $\si$.  If $\pi$ contains no consectutive copies of $\si$ then we say that $\pi$ {\em consecutively avoids} $\si$.  Note that a consecutive copy of $21$ is just a descent while a peak is a consective copy of $132$ or $231$.

Given any finite set of patterns $\Pi$ and  a finite set of positive integers $I$  define
$$
\Pi(I;n)=\{\pi\in\fS_n \mid \text{$\pi$ has a consecutive copy of some $\si\in\Pi$ precisely at the indices in $I$}\}.
$$
Also define the function
$$
\av_\Pi(n) =\#\Pi(\emp;n),
$$
the number of permutations in $\fS_n$ consecutively avoiding all permutations in $\Pi$.  
Given $\Pi\sbe\fS_k$  say that $\Pi$ is {\em nonoverlapping} if for any (not necessarily distinct) $\si,\tau\in\Pi$ and any $l$ with $1<l<k$ the prefix of $\si$ of length $l$ is not order isomorphic to the suffix of $\tau$ of length $l$.
We will now prove our analogue of Theorems~\ref{th:mac} and~\ref{th:bbs} in this setting.
\bth
\label{th:consec}
Let $\Pi\sbe\fS_k$ be a nonoverlapping set of patterns and let $I$ be a finite set of positive integers.  Then for all $n\ge m+k-1$ we have
$\#\Pi(I;n) \in V_\Pi$ where $V_\Pi$ is the vector space of all $\bbQ$-linear combinations of functions in the set
$$
\{n^k \av_\Pi(n+l) \mid k\in\bbZ_{\ge0},\ l\in\bbZ\}.
$$
\eth
\bprf
We induct on $m$.  We have  $\#\Pi(\emp;n)=\av_\Pi(n)$ and so the result clearly holds when $m=0$.
For $m\ge1$,  consider the set $P$ of permutations $\pi\in\fS_n$ which can be written as a concatenation $\pi=\pi'\pi''$ such that $\pi'\in\Pi(I^-;m)$ and $\pi''\in\Pi(\emp;n-m)$.  
Since $\Pi$ is nonoverlapping,  copies of consecutive patterns from $\Pi$ in $\pi$ occur  at  the positions in $I^-$ and possibly also at exactly one of the indices $m,m-1,\dots,m-k+2$.  It follows that 
$$
\#P=\#\Pi(I^-;n)+\#\Pi(I;n)+\sum_{i=1}^{k-2} \#\Pi(I^-\cup \{m-i\};n).
$$

We can also construct the elements of $P$ as follows.  Pick the $m$ elements of $[n]$ to be in $\pi'$ which can be done in $\binom{n}{m}$ ways.  Arrange those elements to have consecutive copies of elements of $\Pi$ at the indices of $I^-$ which can be done in $\#\Pi(I^-;m)$ ways.  Finally, put the remaining elements in $\pi''$ so that it avoids consecutive copies of elements of $\Pi$ which can be done in $\av_\Pi(n-m)$ ways.  Equating  the two counts for $P$ and rearranging terms we get
$$
\#\Pi(I;n)=\binom{n}{m}\av_\Pi(n-m) \#\Pi(I^-;m)-\#\Pi(I^-;n)-\sum_{i=1}^{k-2} \#\Pi(I^-\cup \{m-i\};n),
$$
from which the theorem follows by induction.
\eprf

Note that if $\Pi=\{21\}$ then $\av_\Pi(n)=1$ for all $n$.  So $V_\Pi=\bbQ[n]$ and thus Theorem~\ref{th:mac} is a special case of the previous result.  On the other hand, if $\Pi=\{132,231\}$ then $\av_\Pi(n)=2^{n-1}$ which explains the appearance of the power of $2$ in Theorem~\ref{th:bbs}.  Theorem~\ref{th:consec} suggests that there might be other sets of patterns which would yield interesting enumerative results, and that such sets could be found by looking at $\Pi$ such that the numbers $\av_\Pi(n)$ have nice combinatorial properties.  

\medskip

{\bf (2) The sequence $(a_k(I))$.}  On reading a version of this paper on the arXiv, Ferenc Bencs~\cite{ben:scs}  has found a proof of Conjecture~\ref{con:lc}.  But there is a stronger condition which could also  be investigated.  Consider a finite, real sequence $(a_k)_{0\le k\le n}$ and the corresponding generating function $f(x)=\sum_{k\ge0} a_k x^k$.  It is well known that if the $a_k$ are positive and $f(x)$ has only real roots then the original sequence is log-concave.  However, if one takes $I=\{1,3\}$ then the corresponding generating function is $f(x)=2x^3+6x^2+5x$ which has complex roots.  So this stronger condition does not always apply to the $(a_k(I))$ sequence

\medskip

{\bf (3) Remarks on Conjecture~\ref{-1 basis}.}  Bencs~\cite{ben:scs} has proved this conjecture as well.  His argument is inductive, using the recursions we derived in Section~\ref{sec:rec} as well as Proposition~\ref{prop:c_0} as the base case.   It would be very interesting to prove nonnegativity by finding a combinatorial interpretation of the $c_k(I)$.  Also, one can now further improve the bounds of the roots of $d(I;n)$ on the left side of the $i$-axis by using the linear algebraic method from Section~\ref{sec:roo} on the binomial basis centered at $-1$.

\medskip

{\bf  (4) Limiting behavior of roots.} Bencs~\cite{ben:scs} has proved  a result about the behavior of the roots of $d(I;n)$ for certain sets $I$.  Given $I$, consider the set  $I^k=I\cup\{m+1,m+2,\dots,m+k\}$.
Using Neumaier's Gershgorin-type results on location of polynomial roots~\cite{neu:ecz}, Bencs has demonstrated the following.
\bth
Suppose I is a finite set of positive integers with $m-1\not\in  I$.  Then as $k\ra\infty$ the roots of $d(I^k;n)$ converge to 
$[0,m+k]-\{m-1\}$.\hqed
\eth
\medskip

{\em Acknowledgement.}  We wish to thank Marcelo Aguiar for asking the question that lead to this research.  Thanks also to Marcelo Aguiar, Jason Brown, Petter Br\"and\'en, Ira Gessel, John Stembridge, and Richard Stanley for helpful discussions and  useful references. {A. Diaz-Lopez thanks the AMS and Simons Foundation for support under the AMS-Simons Travel Grant.} {P.~E.~Harris was partially supported by NSF grant DMS--1620202.} {M. Omar thanks the Harvey Mudd College Faculty Research, Scholarship, and Creative Works Award.} 



\newcommand{\etalchar}[1]{$^{#1}$}

\end{document}